\documentclass[11pt]{amsart}
\usepackage{amsmath}
\usepackage{amssymb}
\newtheorem{theorem}{\bf Theorem}[section]

\newtheorem{proposition}{\bf Proposition}[section]

\newtheorem{remark}{Remark}[section]

\newcommand{\ad}{\mbox{ad }}
\newcommand{\Ad}{\mbox{Ad }}

\newcommand{\diag}{\mbox{\bf diag} }
\begin{document}

\title{GLOBAL CONTROLLABILITY OF MULTIDIMENSIONAL RIGID BODY BY FEW TORQUES}

\author{A. V. Sarychev$^*$}

\address{Dipartimento di Matematica per le Decisioni, Universit\`a di Firenze,\\
v.~C.Lombroso 6/17, Firenze, 50134, Italy\\
$^*$E-mail: asarychev@unifi.it}

\begin{abstract}
We study global controllability of 'rotating' multidimensional rigid body (MRB) controlled by application of few torques.
Study  by methods of geometric control requires analysis of algebraic structure
introduced by the quadratic term of Euler-Frahm equation. We discuss  problems,
which arise in the course of this analysis, and
establish several global controllability criteria for damped and non damped cases.
\end{abstract}


\maketitle

\section{Introduction}\label{intr}
In recent work \cite{AS1,AS2,Ro} one studied  controllability of
Navier-Stokes (NS) equation,  controlled by forcing applied to few modes on a 2D domain.
Geometric control approach has been employed for establishing approximate controllability criteria
   for NS/Euler  equation on 2D torus, sphere, hemisphere, rectangle and  generic Riemannian surface with boundary.

In the present contribution we address  controllability issues for a finite-dimensional "kin"  of NS equation - Euler-Frahm equation for rotation of multidimensional rigid body ({\bf MRB}) subject to few controlling torques and to possible damping. The equation  evolves on $so(n)$. We formulate global controllability criteria which are structurally stable with respect to
the choice of 'directions' of controlled torques.

According to geometric approach to studying controllability,   one starts
with a system controlled by low-dimensional input and
proceeds with a sequence of {\em Lie extensions}  (\cite{Jur,Sa06}) which add to the system new controlled
vector fields. The latter are calculated via iterated Lie-Poisson
brackets  of  the controlled vector fields and the drift (zero control vector field).
The core of the method and the main difficulty is in finding proper Lie extensions and in tracing results of their implementation.

The Lie extension employed in \cite{AS1,AS2,Ro} for  studying controllability of  NS equation, and similar one used  equation below (see Subsection~\ref{41}), involves  double Lie bracket of drift  vector field with a couple of {\em constant} controlled vector fields (they are identified with their values,  or {\em directions} belonging to $so(n)$).
At least one of the directions must be a {\em steady state} of {\bf MRB}, i.e. an  'equilibrium points' of  Euler-Frahm equation.

 The double Lie bracket results in constant controlled vector field (extending direction); the  correspondence between couple of original controlled directions and the extending one defines bilinear operator $\beta$ on $so(n)$. More extending controlled directions are obtained by {\em iterated} application of $\beta$.
     For proving {\em global controllability of {\bf MRB}} we must verify {\em saturating property} - coincidence of the set of  extending directions with $so(n)$ after a  number of  iterations.

Tracing the iterations is by no means easy.   For NS equation all  cases, successfully analyzed in  \cite{AS1,AS2,Ro}, are related to an explicit description of the basis of steady states and to specific representation of the
  operator $\beta$  with respect to this basis.
  The results,  so obtained, are heavily dependent on choice of original controlled directions and on geometry
   of the domain where the NS equation evolves.

Below we manage to establish several controllability criteria for damped and non damped {\bf MRB} controlled by one, two or three torques.  We pay special attention to deriving criteria which are structurally stable with respect to perturbation of (some of) the controlled directions.

\section{Euler Equation for Generalized Rigid Body and Euler-Frahm equation for MRB}
\label{frahm}
We follow \cite{Arn} for definition of 'generalized rigid body'. Let $\mathcal{G}$ be a Lie group, $\mathfrak{g}$ its Lie algebra  and let left-invariant Riemannian metric on $\mathcal{G}$ be defined by scalar product $\langle \cdot , \cdot \rangle$ on $\mathfrak{g}$.

Introduce $\mathcal{I}: \mathfrak{g} \mapsto \mathfrak{g}^*$ -  a symmetric operator, which corresponds to the Riemannian metrics by  formula: $\langle \xi , \eta \rangle=\mathcal{I} \xi | \eta$, where $\cdot | \cdot $ is the natural pairing between $\mathfrak{g}$ and $\mathfrak{g}^*$.
The operator $\mathcal{I}$ is called {\em inertia operator} of generalized rigid body.

The trajectory of the motion of generalized rigid body is a curve $g(t) \in \mathcal{G}$. {\em Angular velocity}, corresponding to this motion is:
$  \Omega=L_{g^{-1}*}\dot{g} \in  \mathfrak{g}$,
 where $L_g$ is left translation by $g$.
The image of angular velocity $\Omega$ under $\mathcal{I}$ is {\em angular momentum} $M \in \mathfrak{g}^*$. Energy of the body equals $\langle \Omega , \Omega \rangle =M |\Omega$.

{\em Euler equation} for the motion of generalized rigid body is
$ \dot{\Omega}=\mathcal{B}(\Omega, \Omega),$
where bilinear operator $\mathcal{B}: \mathfrak{g} \times \mathfrak{g} \mapsto \mathfrak{g}$ is defined by formula:

\begin{equation}\label{defb}
    \langle  [a,b], c \rangle = \langle \mathcal{B} (c,a),b \rangle ,
\end{equation}
$[\cdot , \cdot ]$ staying for Lie-Poisson bracket in $\mathfrak{g}$.

    {\bf MRB} is particular case of generalized rigid body, where the Lie group $\mathcal{G}=SO(n)$, and  {\em angular velocities} $\Omega \in \mathfrak{g}=so(n)$ are skew-symmetric matrices.

Identifying $so(n)$ with $so^*(n)$ by means of Killing form,  we may think of  momentum $M$
as of skew-symmetric matrix. Then the inertia operator is a map
$
\mathcal{I}_C: \Omega   \mapsto  \left(\Omega  C + C \Omega \right)=
M \in
so(n) \stackrel{\mathcal{K}}{\cong} so^*(n), $
where $C$ is some positive semidefinite matrix.

 Operator $\mathcal{I}_C$ is symmetric with respect to Killing form and
  is invertible  (Sylvester theorem),
whenever  $C$ is positive definite.

We compute  $\mathcal{B}$ according to
(\ref{defb}) ($[\cdot, \cdot]$ being matrix commutator):
\[ \mathcal{B}(\Omega^1,\Omega^2)=
 \mathcal{I}^{-1}_C\left[\mathcal{I}_C\Omega^1, \Omega^2 \right], \ \mathcal{B}(\Omega,\Omega)=
 \mathcal{I}^{-1}_C\left[\mathcal{I}_C\Omega, \Omega \right]=\mathcal{I}_C^{-1}[C,\Omega^2].\]
{\em Euler-Frahm equation} for the motion of {\bf MRB} is:
\begin{equation}\label{efad0}
\dot{\Omega}= \mathcal{I}^{-1}_C\left[\mathcal{I}_C\Omega, \Omega \right]=\mathcal{I}_C^{-1}[C,\Omega^2],
\end{equation}
The motion, subject to  damping, is described by the equation
\[\dot{\Omega}=\mathcal{I}_C^{-1}[C,\Omega^2] - \nu \Omega , \ \nu \geq 0.\]

\section{Controllability of rotating MRB: problem setting and main results}
{\em Controlled} rotation of {\bf MRB} is described by equation
 \begin{equation}\label{crot}
\dot{\Omega}=\mathcal{I}_C^{-1}[C,\Omega^2] - \nu \Omega +\sum_{i=1}^rG^i
u_i(t), \ \nu \geq 0, \ G^i \in so(n).
 \end{equation}

We are interested in {\em global controllability} of (\ref{crot}),  meaning
that for any  $\tilde{\Omega},\hat{\Omega} \in so(n)$ system (\ref{crot}) can be steered  from $\tilde{\Omega}$ to $\hat{\Omega}$ in some time $T \geq 0$.
We are interested in achieving global controllability  by small number of controls; we
prove that $r$ can be taken $\leq 3$ for all $n \geq 3$.

Equation (\ref{crot}) is particular case of control-affine system with quadratic(+linear) {\em drift vector field} and {\em constant controlled vector fields}.

    The following {\em genericity condition} is assumed to hold furtheron: symmetric  matrix $C$   {\em is positive definite and has distinct  eigenvalues}.

Our first result claims {\em global controllability of {\bf MRB} by means of two controlled torques}.

\begin{theorem}
\label{maida}
There exists a pair of directions $G^1,G^2 \in so(n)$ (depending on $C$),  such that the system (\ref{crot})
with $r=2$ is globally controllable. $\Box$
\end{theorem}

 The proof of this Theorem, sketched below, is based on direct computation of Lie extensions in specially selected basis, related to $C$. More difficult is formulating  criteria, which are {\em structurally stable} with respect to perturbation of controlled directions.

 We start with {\em non damped} {\bf MRB}, controlled by one torque.
 In this case -  given  {\em recurrence} of  Euler-Frahm dynamics (\ref{efad0}) -
 {\em bracket generating property} suffices for guaranteeing global controllability.
This property means that evaluations (at each point) of iterated Lie brackets of drift and controlled vector fields
span $so(n)$.
Given high dimension of $so(n)$,
verification of the bracket generating property for {\em generic}  controlled direction is nontrivial task.
We do this analyzing linearization of quadratic Euler operator.
The result is

\begin{theorem}
  \label{undcon}
  For generic $G \in so(n)$ the system
  $\dot{\Omega}=\mathcal{I}_C^{-1}[C,\Omega^2] +G u(t),$
  is globally controllable, also if control is bounded: $|u| \leq b, \ b>0. \ \Box$
\end{theorem}

We now pass to the damped case. Our  method requires one of the controlled directions to
    be {\em steady state} for  {\bf MRB}.
    Recall that steady state or steady direction  of  {\bf MRB}
    is equilibrium point of (\ref{efad0}) - a matrix  $\hat{G}$ for which
        $[\mathcal{I}_C \hat{G},\hat{G}]=[ C ,\hat{G}^2]=0 $. Matrix $\hat{G}$ is {\em principal axis} of {\bf MRB}, if $\mathcal{I}_C \hat{G}=\mu  \hat{G}, \ \mu \in \mathbb{R}$.  These two sets coincide for $n=3$, while for
         $n \geq 4$ the
         set of steady directions is much richer.

The results obtained for the damped case differ for {\em odd} and {\em even} $n$.

\begin{theorem}
\label{dampodd}
Let $r=2, \ n$ be odd  in (\ref{crot}). For generic stationary direction $G^1$ and generic $G^2 \in so(n)$
 the system (\ref{crot}) is globally  controllable. $\Box$
\end{theorem}

An additional symmetry in the case of even $n$, obliges one to involve additional controlled direction
for achieving global controllability.

\begin{theorem}
\label{dampeven}
Let $r=3, \ n$ be even in (\ref{crot}). For generic  stationary direction $G^1 \in so(n)$  and  generic pair  $(G^2,G^3)$ of directions
  the system (\ref{crot}) is globally controllable. $\Box$
\end{theorem}

Generic element of a subset $W \subseteq so(n)$ means an element of  open dense subset of $W$ in induced topology.

\section{Sketch of the proof of Theorem~\ref{maida}}

\subsection{Key Lie extension}\label{41} {\em Lie extensions}  mean finding vector fields
$X$, which  are  {\em compatible} with control system, in the sense that
 closures of attainable
sets of the control system are invariant for $X$.
 If one is able to prove global controllability
 of the system  {\em extended} by some compatible vector fields,
 then controllability of the original system can be
concluded by  standard argument.

Key Lie extension, we employ, is described by the following
\begin{proposition}\label{lext}
Let  for control system
  \begin{equation}\label{f012}
    \dot{x}=f(x)+\tilde{g}(x)u+\bar{g}(x)v,
\end{equation}
evolving on a manifold $Q$,  hold the relations
\begin{equation}\label{skob}
 \{\tilde{g},\bar{g}\}=0, \ \{\tilde{g},\{\tilde{g},f\}\} =0,
\end{equation}
($\{\cdot, \cdot\}$ stays for Lie brackets of vector fields on $Q$). Then the system
$\dot{x}=f(x)+\tilde{g}(x)u+\bar{g}(x)v+\{\bar{g},\{\tilde{g},f\}\}(x)w$
is  Lie extension of (\ref{f012}). $\Box$
\end{proposition}

\begin{remark}
Vector fields $\pm \{\bar{g},\{\tilde{g},f\}\}$ are {\em extending controlled vector fields}; they are also compatible with (\ref{f012}).  $\Box$
\end{remark}

We will repeatedly employ Proposition~\ref{lext} for extending control system (\ref{crot}).
At each step the first of the relations (\ref{skob})
will be  trivially satisfied since all original and extending controlled vector fields will be
constant. For drift vector field
$f(\Omega)=\mathcal{I}_C^{-1}[C,\Omega^2]$ in (\ref{crot}), and constant controlled vector field
$\tilde{g} \equiv \tilde{G} \in so(n)$ , the Lie bracket $\{\tilde{g},\{\tilde{g},f\}\} \equiv \mathcal{I}_C^{-1}[C,\tilde{G}^2]$ is constant vector field. The second relation
(\ref{skob}) would hold if and only if
$\tilde{G}$ is  steady state. When repeating the extension it is crucial to guarantee
at each step disponibility of steady state controlled direction.

        For  two constant controlled vector fields $\tilde{g} \equiv \tilde{G},\bar{g} \equiv \bar{G}, \ \tilde{G},\bar{G} \in so(n)$ the value of constant {\em extending controlled vector field} $\{\bar{g},\{\tilde{g},f\}\}$ is
  \begin{equation}\label{ggb}
\beta(\tilde{G},\bar{G})=\mathcal{I}^{-1}[C,\tilde{G}\bar{G}+\bar{G}\tilde{G}];
\end{equation}
 (\ref{ggb}) defines symmetric bilinear operator $\beta$ on $so(n)$.

\subsection{Algebra of principal axes and controllability proof}

 Diagonalize matrix  $C$ presenting it as $C=\Ad S D=SDS^{-1}$ with $S$ orthogonal and  $D=\diag \{I_1, \ldots, I_n\}, \ I_1 < I_2 <
\cdots <I_n$.

Introduce matrices $\Theta^{rs}={\mathbf 1}_{rs}-{\mathbf 1}_{sr} \in so(n), \ (1 \leq r < s \leq n),$ with
${\mathbf 1}_{rs}$ being  matrix with (the only nonvanishing) unit
element at row $r$ and column $s$.
Matrices $\Omega^{rs}=\Ad S \Theta^{rs}$ turn out to be 'eigenvectors' of the operators $(\ad C)$ and $\mathcal{I}_C$. They form set of principal axes of the {\bf MRB}.

'Multiplication table' for  $\beta$ with respect to the basis $\Omega^{rs}$ is
\begin{equation*}
\beta(\Omega^{rs},\Omega^{rs})=0, \ \beta(\Omega^{rs},\Omega^{r\ell})=(I_\ell-I_s)(I_s+I_\ell)^{-1}\Omega^{s\ell},
\end{equation*}
$\beta(\Omega^{rs},\Omega^{k\ell})=0$, whenever  $r,s,k,\ell$ are
distinct.

Take \[G^1= \Omega^{12} \ \mbox{- principal axis}, \
G^2=\Omega^{23} + \Omega^{34}
+\cdots +\Omega^{n-1,n}.\]

     It suffices to prove that iterated applications of $\beta$ to  $G^1,G^2$  result in a basis of $so(n)$,
because then the extended system would possess full-dimensional
input and therefore would be globally controllable. The original system (\ref{crot}) would be  globally controllable as well.

According to the multiplication table
$G^3=\beta(G^2,G^1)=\beta(\Omega^{12},\Omega^{23})$
 coincides up to a multiplier with principal axis $\Omega^{13}$.
Calculating subsequently extending controlled 'directions'
$G^{i}=\beta(G^{i-1},G^2), \ i >
2,$ we see that all $G^{i}$ coincide up to a nonzero
multiplier with $\Omega^{1,i}$, i.e.  are principal axes.
Also $\beta(\Omega^{1i},\Omega^{1k})$ coincides up to a multiplier with $\Omega^{ik}$; this means
that iterating applications of $\beta$ to $G^1,G^2$ generate basis of $so(n)$.

\section{Concluding remarks}

1. As one can see  proof of Theorem~\ref{maida} is "rigid construction",  based on specific  choice of controlled directions  and on computation of iterated Lie extensions with respect to specific basis of principal axes of {\bf MRB}. If one perturbs  one of the original controlled directions the  constructions fails, as far as first Lie extension does not result in  new stationary direction of {\bf MRB}, and the  Proposition~\ref{lext} can not be iterated.

This rigidity of controllability criteria with respect to the choice of controlled directions,  manifested itself also in previous study of approximate controllability of NS system on particular 2D domains  (\cite{AS1,AS2}). It does not seem natural, and is rather related to the proposed method.

Indeed one would expect structural stability of controllability criteria and this is achieved in the formulations of   Theorems~\ref{undcon},\ref{dampodd},\ref{dampeven},   which are structurally  stable with respect to the choice of (some of) the controlled directions. The method for establishing these  criteria  differs from the previous one. It is based on study of linearization at a steady state of Euler operator for {\bf MRB}. The proofs will appear elsewhere.

Besides its  interest for studying controllability of {\bf MRB}, 
the method  can be extended onto infinite-dimensional case, and be applied  to 
controlled NS/Euler equation for fluid dynamics on general 2D and 3D domains. 
The results will appear in further publications.

2. Publication \cite{Der} studies controllability of
 non damped {\bf MRB} by using of a pair of controlled 'flywheels' -  different
 type of "internal-force  controls", with dynamics described by    bilinear control system on Lie group.

\end{document}